\documentclass{elsarticle}
\usepackage{amsfonts} 
\usepackage{graphicx}
\usepackage{float}
\usepackage{amsmath}
\usepackage{url}
\newtheorem{thm}{Theorem}
\newtheorem{lem}{Lemma}
\newtheorem{rem}{Remark}

\newtheorem{exmp}{Example}

\begin{document}
\begin{frontmatter}
\author{Lajos L\'aszl\'o}
\address{Department of Numerical Analysis, 
E\"otv\"os L\'or\'and University, Budapest}
\title{On the Grace-Danielsson inequality for tetrahedra}

\begin{abstract} 
The difference between the (squared) sides of the Grace-Danielsson
inequality for tetrahedra will be represented as a sum of 
two nonnegative terms.
This gives another proof of the inequality. Examining the 
denominator allows us to characterize the infinite triangular 
prism as a degenerate tetrahedron. 
We give conditions for equality (for a zero gap) as well, 
and some examples are included.
\end{abstract}

\begin{keyword} tetrahedron \sep Grace-Danielsson inequality 
\sep inradius \sep circumradius 
\MSC{51M04, 51M16}
\end{keyword} \end{frontmatter}

\section{Introduction}

  It is a classic result of 18th century mathematics (Chapple, Euler)
that
   \begin{equation} \label{Eul} d^2=R(R-2r) \end{equation}
holds for the distance $d$ between the circumcenter and the incenter 
of a triangle with circumradius $R$ and inradius $r.$  
Although a wish to generalize it to tetrahedra failed (cf. Gergonne, 
Durrande, 19th century), an inequality
   \begin{equation} \label{GrD} d^2\le(R+r)(R-3r) \end{equation}
is still valid for all tetrahedra with circumradius $R,$ inradius $r$
and distance $d$ between them (see Grace \cite{Grace}, 
Milne \cite{Milne} for further references. Milne and others cite
Danielsson \cite{Daniel} which is not readily available).
The  ideas of Milne, using quantum information theory can now
be found in his recently published dissertation 
\cite{MilneDissert} as well. 

Now, it is a natural question, how big is the difference between 
the right and left hand sides. Our aim is to become a
representation, from which inequality (\ref{GrD}) 
evidently follows.
To this we first rewrite the inequality as
   \[ R^2-d^2-3r^2 \ge 2 r R, \]
and square both sides to avoid irrationality, getting thus a 
purely algebraic form. 

When we  met this inequality, we tried to prove it by means of
semidefinite programming, more concretely, by sum of squares
programming. However, this attempt failed due to perpetual
``Out of memory'' problems. This is why we were forced 
to have a closer look at this representation problem and
to write our own routines in Maple and Matlab
(the latter for quickening and/or checking purposes).
Thus the result obtained also serves as an illustration for the 
positive answer for Hilbert's 17th problem \cite{Hilb}.

As regards variables used, we prefer to choose the coordinates 
of the three vertices (forming the so-called basic face) 
and the tangent point of the insphere lying on the basic face,
as well as the inradius $r$ -- and calculate the coordinates of 
the remaining vertex, the circumradius $R$ and the distance $d.$
With this choice we get the wanted rational representation.
Our method can be hence considered -- apart from the use of 
programming languages -- elementary.

Using coordinates, as is known, goes hand in hand with long 
calculations, however the relevant formulas here can be
well managed by means of the symbolic programming
language Maple. 
The details will be given in the proof of Theorem 1, 
followed by some special cases
(concerning the choice of the tangent point),
illustrative examples, and by investigating the degeneracy also
in two dimensions (for a triangle instead of for a tetrahedron).

\section{The main theorem} 

At first we formulate the theorem, giving a two-term
representation for the gap, where the quantities on the right 
hand side will be explained in the course of the proof.

\begin{thm}
 For a tetrahedron with inradius $r,$ circumradius $R,$ and
distance $d$ between the incenter and circumcenter we have
  \begin{equation} \label{th} 
  \big(R^2-d^2-3r^2\big)^2 - (2 r R)^2 = r^2 \frac{ 
  (u_1 r^2+v_1)^2 + (u_2 r^2+v_2)^2} {a_0 (A-B\, r^2)}  
\end{equation}
with polynomials $u_1, v_1, u_2, v_2, A, B, a_0,$ 
where $A, B, a_0$ are positive, and all these quantities 
depend only on one face, called the basic triangle. 
\end{thm}
{\bf Proof.} 
First we describe our method in a more detailed form. Let 
$x, y, z$ be the vertices of the {\it basic} triangle, considered lying 
on the horizontal plane. 
Let $c$ be an interior point of it, and $r>0$ be a given number.  
We will find the fourth vertex $w$ such that the insphere of the 
tetrahedron $\{x, y, z, w\}$ has radius $r$ and touches the basic 
triangle at its inner point $c.$

To this aim we need to draw the three tangent planes and then to find 
their common point $w.$ Since a too large value $r$ contradicts the
requirements, it is essential to know the largest possible, 
the {\it critical} value of the inradius.
In case of the critical situation all the three edges are parallel, 
forming
a semi-infinite triangular prism. (Example \ref{crit} illustrates this 
behavior, while the planar version is discussed in Lemma \ref{planar}.)

 Denote by $x= (x_1, x_2, 0), \, y=(y_1, y_2, 0), \,z=(z_1, z_2, 0)$ 
the vertices of the basic triangle (oriented counterclockwise), by 
$(c_1, c_2, 0)$ the interior point chosen, and let $r>0$ be a given number.
(Note that we write ($x_1, x_2, x_3$) for coordinates, instead of 
 ($x, y, z),$ used in some cases.) 

As a first step, we determine the three tangent points  $X, Y, Z$ 
on the insphere. 
Here $X, Y, Z$ are opposite to vertices $x, y, z,$
resp. Calculate then the intersection point $w$ of the tangent 
planes spanned by the triangles
$\{x, y, Z\}, \ \{y, z, X\} $ and $\{z, x, Y \}.$ 
Although the first two coordinates $w_1, w_2$ of $w$ are complicated, 
the third -- and most relevant -- can be handled well.
With the quantities  
\begin{eqnarray*} 
 a_0&=& 2\, \mathrm{area}\, \Delta xyz = 
 x_1 y_2+y_1 z_2+z_1 x_2 -y_1 x_2-z_1 y_2-x_1 z_2, \\
 a_x&=& 2\, \mathrm{area}\, \Delta cyz= 
 c_1 y_2\, +y_1 z_2\, +z_1 c_2  -y_1 c_2 -z_1 y_2 -c_1 z_2, \\
 a_y&=& 2\, \mathrm{area}\, \Delta xcz= 
 x_1 c_2+c_1 z_2+z_1 x_2 -c_1 x_2-z_1 c_2-x_1 z_2, \\
 a_z&=& 2\, \mathrm{area}\, \Delta xyc= 
 x_1 y_2+y_1 c_2+c_1 x_2 -y_1   x_2-c_1 y_2-x_1 c_2,
 \end{eqnarray*}
and 
\[A=a_x a_y a_z, \quad 
B= \|x\|^2 a_x+ \|y\|^2 a_y+ \|z\|^2 a_z- \|c\|^2 a_0\]
we get the formula 

\begin{equation} \label{w3} w_3=\frac{2\, r A}{A-B r^2}. 
\end{equation}

Here the quantity $A$ is -- as a product of three triangle areas -- 
obviously positive, while the same property for $B$ will be 
proved in Lemma 1 below. Then it is seen that $w_3>0$ for $r$ sufficiently small, and that its critical value is 
\begin{equation*}  r_{crit} = \sqrt{A/B},\end{equation*}
where the tetrahedron becomes a prism. Notice that by means of
this critical value the third coordinate of $w$ can be rewritten as
\begin{equation} \label{w3byrcrit} w_3=\frac{2\,r\, r_{crit}^2}
{r_{crit}^2-r^2}. 
\end{equation}
Now we determine the circumcenter $o.$ For this case -- in 
contrast with vertex $w$ -- the first and second coordinates are 
relatively simple, 
while the third one is difficult (but not needed here). 
We have
\[ 2\, a_0\, o_1=\|x\|^2 (y_2-z_2)+\|y\|^2 (z_2-x_2)+\|z\|^2 
(x_2-y_2), \]
\[ 2\, a_0\, o_2=\|x\|^2 (z_1-y_1)+\|y\|^2 (x_1-z_1)+\|z\|^2 
(y_1-x_1). \]
Observe that the orthogonal projection $(o_1, o_2, 0)$  of the circumcenter 
$(o_1, o_2, o_3)$ of the tetrahedron coincides with
the circumcenter of the basic triangle, due to elementary 
considerations, hence $o_1, o_2$ do not depend on $w.$

Now we are in the position to calculate the circumradius $R=\|o-x\|$ 
and the distance $d=\|o-c\|$ of the circumcenter and the incenter.
In possession of these, a quite circumstantial manipulation in Maple 
is needed to get a more concise form for the gap, resulting in the two 
term-representation of the theorem. 
Of the remaining four variables, $u_1$ and $u_2$ have a 
fairly simple form:
\begin{equation} \label{ugyan} u_1=4 a_0 (c_1-o_1), 
\quad  u_2=4 a_0 (c_2-o_2), \end{equation}
while $v_1$ and $v_2$ are polynomials of degree five with $108-108$ 
terms, but they can be rewritten into a something shorter 
sum-of-products form, see Appendix.
The subtask of determining $v_1$ and $v_2$ knowing the $u_i$-s
leads to the following.
Assume that $\alpha=u_1^2+u_2^2, \ 
\beta, \ \gamma$ are known and $v_1, \, v_2$ are asked to satisfy
\[ \alpha r^4 + \beta r^2 + \gamma = 
 (u_1 r^2+v_1)^2 + (u_2 r^2+v_2)^2, \]
then the solution is given by
\[v_1=\frac{u_1 \beta \pm u_2 \mathrm{dis}}{2  \alpha}, \ \
 v_2=\frac{u_2 \beta \mp u_1 \mathrm{dis}}{2  \alpha}, \quad  
 \mathrm{dis}=\sqrt{4 \alpha \gamma -\beta^2}. \] 

Since the discriminant dis (having originally as much as
370900 terms!) proves to be a complete square, 
we arrive at the representation (\ref{th}), expressing the gap
as a rational function of the variables. 
The theorem is proved. $\square$ \medskip

Since the right hand side in  (\ref{th}) is non-negative, we have 
thus another proof for the Grace-Danielsson inequality.
It remains to prove the positivity of $B.$   \medskip

\begin{lem} \label{trans} The polynomial $B$ in the denominator of 
(\ref{th}) is positive, further, $B$ is independent of translation.
\end{lem}

{\bf Proof.} 
Let 
\[ t_x=\frac{a_x}{a_0}, \quad  t_y=\frac{a_y}{a_0}, \quad  
t_z=\frac{a_z}{a_0}, \]
then these positive numbers sum up to one and we have to show
 \[ \|c\|^2 \le \|x\|^2 t_x+\|y\|^2 t_y+\|z\|^2 t_z. \]
The well known barycentric formula
\[ c=\frac{a_x x+a_y y+a_z z}{a_x + a_y + a_z} =
 t_x x+t_y y+t_z z\]
implies in our case
\[  \|c\|^2 < \big(t_x \|x\| +t_y \|y\| +t_z \|z\| \big)^2, \]
whence the Cauchy-Schwarz inequality applied for vectors
\[ (\sqrt{t_x}, \sqrt{t_y},\sqrt{t_z}), \quad
 (\|x\| \sqrt{t_x},\, \|y\| \sqrt{t_y},\, \|z\| \sqrt{t_z}) \]
yields the positivity of $B.$
To prove the second assertion, translate now the vertices 
by $h$ to get
\begin{eqnarray*} 
B(h)&=& \|x-h\|^2 a_x+\|y-h\|^2 a_y+\|z-h\|^2 a_z -
\|c-h\|^2 a_0 \\
&=& \|x\|^2 a_x+\|y\|^2 a_y+\|z\|^2 a_z -\|c\|^2 a_0 \\
&-& 2 h^T \big(a_x x+a_y y+a_z z - a_0 c \big) \\
&+& \|h\|^2 \big(a_x + a_y + a_z - a_0 \big). \end{eqnarray*}
The coefficient of $\|h\|^2$ is obviously zero by additivity of the 
area, while the linear term disappears because of the basic property 
of barycentric coordinates.
Thus $B(h)\equiv B(0)=B,$ which was to be shown. $\square$ 
\medskip

Now we give a numerical example using a Heronian tetrahedron,
for which the essential quantities all are rational 
(cf. \cite{mathoverflow}), justifying thus
 the use of Maple. 
\begin{exmp}
Let the vertices of the basic triangle, and the tangent point
of the insphere  be
$ x=(0,0,0), \ y =(154, 0, 0), z = (55, 132, 0), \ c = (90, 48, 0), $
and choose $r=10.$  
Then the fourth vertex and the circumcenter become
\[w=\Big(\frac{215490}{2309}, \frac{339416}{6927}, 
\frac{49280}{2309}\Big), \quad o=\Big(77, \frac{363}{8}, 
-\frac{15818598389}{93098880}\Big). \]

Further we have \medskip

\noindent $R^2=319462309835987155321/8667401457254400,$ 

\noindent $\, d^2=282073185661355308921/8667401457254400,$

\noindent $f=198873308525/145467,$

\noindent $a_0=20328, \ a_x=3696, \ a_y=9240, \ a_z=7392,$ 

\noindent $A=252444487680, \ B=158802336,$

\noindent $u_1=1057056, v_1=-7868399616, \ $
               $u_2=213444, v_2=-2363251968.$
\end{exmp}

{\bf Question.} 
Since the tetrahedron, a 3-dimensional simplex, has a two-term
gap given by the right hand side of (\ref{th}), one can put the 
question: how many terms (if any) can represent 
the gap for a simplex in $n>3$ dimensions? \smallskip

See to this John Baez's blog \cite{BaezEgan} citing Greg Egan,  
for the concrete form
\[d^2 \le (R+(n-2)r) (R-nr)\]
of the inequality in $n$ dimensions -- or the equivalent, 
``Pythagorean'' form
\begin{equation} \label{Pyth} d^2 + (n-1)^2 r^2 \le (R-r)^2.
\end{equation}

Back to $n=3,$ the next example shows that equality in
 (\ref{GrD}) (or in (\ref{Pyth}))  can occur for 
non-regular tetrahedra, in contrast with Euler's inequality 
$R\ge 2r,$ where equality is valid only for regular triangles.

\begin{exmp} Let the vertices of the basic triangle, 
the inner point chosen, and the inradius be
\[ x=(-1, 0, 0), \ y=(1, 0, 0), \ z=(0, \sqrt{3}, 0), \ 
c=(0, \tfrac{1}{\sqrt{3}}, 0), \ r=\tfrac{1}{2}. \]
From these data the method gives the fourth vertex $w,$ 
the circumcenter $o,$ 
\[ w=\Big(0,\frac{1}{\sqrt{3}},4\Big), \quad  
o=\Big(0,\frac{1}{\sqrt{3}}, \frac{11}{6}\Big), \]
and the further parameters
\[ a_x=a_y=a_z=\frac{2\sqrt{3}}{3}, \ A=\frac{8\sqrt{3}}{9}, 
\ B=\frac{8\sqrt{3}}{3}, \ R=\frac{13}{6}, \ d=\frac{4}{3}. \]
Therefore  (\ref{Pyth}) turns into equality thanks 
to the Pythagorean identity  $3^2+4^2=5^2.$ 
This result also follows by Lemma \ref{lem-uv} below.
\end{exmp}

\begin{section} {Some special cases and examples} \end{section}

The next lemma describes the gap with disappearing $u_i$-s, 
and $v_i$-s, resp. \medskip

\begin{lem} \label{lem-uv}

(a) For a basic triangle touched by the insphere at its circumcenter 
$c_1=o_1, \ c_2=o_2$ the $u_i$s vanish and the gap (\ref{th}) is 
given by 
\begin{equation} \label{eqa} 
\frac{r^2}{64\,a_0^5\,(A-Br^2)}\, \|x-y\|^4 \|y-z\|^4
 \|z-x\|^4 (g_1^2+g_2^2) 
\end{equation}  
with
\begin{eqnarray*} 
g_1=&(x_1^2+3 x_2^2)(z_2-y_2)+(y_1^2+3 y_2^2)(x_2-z_2)+
  (z_1^2+3 z_2^2)(y_2-z_2) \\ &+2 x_1 x_2 (z_1-y_1) + 
  2 y_1 y_2  (x_1-z_1) + 2 z_1 z_2 (y_1-x_1), \\
g_2=&(x_2^2+3 x_1^2)(y_1-z_1)+(y_2^2+3 y_1^2)(z_1-x_1)+
  (z_2^2+3 z_1^2)(x_1-y_1) \\ &+2 x_1 x_2 (y_2-z_2) + 
  2 y_1 y_2 (z_2-x_2) + 2 z_1 z_2 (x_2-y_2). \end{eqnarray*} 

(b) For a basic triangle touched by the insphere at its 
incenter, the $v_i$s vanish and the gap is given by
\begin{equation} \label{eqb}  \frac{16\, a_0\, r^6}{A-B r^2}\, 
\big( (c_1-o_1)^2+(c_2-o_2)^2 \big). \end{equation} 

(c) Furthermore, if the basic triangle is regular with 
$c_1=o_1, c_2=o_2$, then the gap is 0.
\end{lem}  

{\bf Proof.}

(a) The statement follows by (\ref{ugyan}).
 Factorization by Maple gives
\[ v_i=\frac{1}{8\, a_0^2} \, \|x-y\|^2\, \|y-z\|^2\, \|z-x\|^2\, 
g_i, \, \ i=1,2, \]
with the third degree polynomials $g_1, g_2$ above.

(b) If we calculate the incenter (by using barycentric coordinates), 
it appears that $v_1 = v_2 = 0$ holds, and the result follows from
(\ref{th}).

(c) In an equilateral triangle circumcenter and incenter are coincident 
at the center of rotational symmetry, so if the basic triangle is 
equilateral, and the touching point is the center of symmetry, the gap 
is given either by case (b), equation (\ref{eqb}), i.e. is 0 
($c_i = o_i$) 
or by case (a), equation (\ref{eqa}). Both equations having to give 
the same result, it implies that the term ($g_1^2 + g_2^2$) in  
(\ref{eqa})
has to be 0 because no other terms of this equation can cancel. 
Therefore the polynomials $g_1$ and $g_2$ of (3) are nil.
This can be checked using Maple: to force the basic triangle 
to be regular, we further substitute
\[ z_1= \tfrac{1}{2} \big (x_1+y_1+\sqrt{3}\,(y_2-x_2)\big), 
\ \ z_2= \tfrac{1}{2} \big(x_2+y_2+\sqrt{3}\,(x_1-y_1)\big) \]
in $g_1, g_2,$ to get $g_1=g_2=0.$ 
  \ $\square$

 \begin{figure}[h!]
 \includegraphics[scale=0.65]{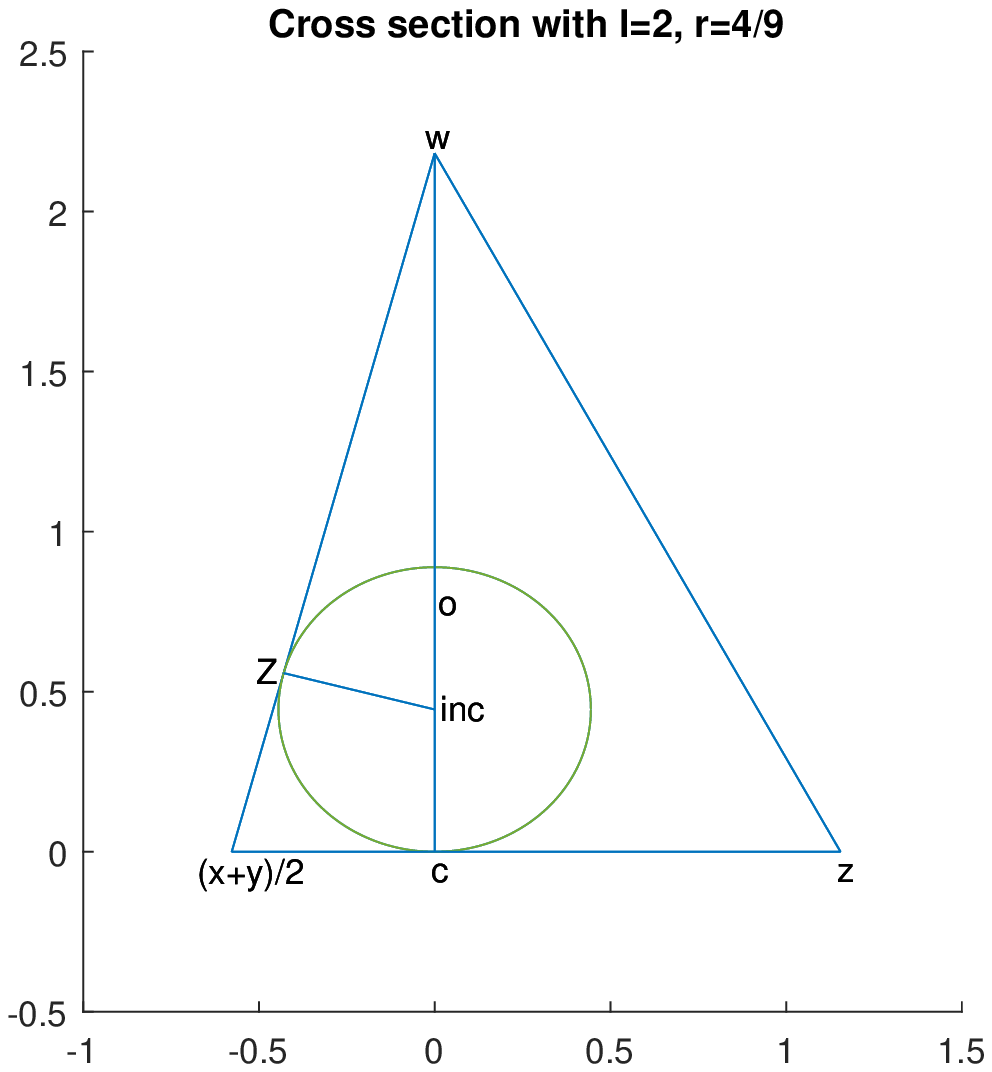}
 \centering
 \centerline{Figure 1. The circumcenter is inside the tetrahedron}
 \end{figure}
 
\begin{rem}
It was the referee's idea to give a Mapleless proof for the zero gap property (c). Also, he provided a proof (essentially part (c2) 
below), where he obtained (\ref{wnew}) below 
by using (\ref{w3byrcrit}), as a consequence of (\ref{w3}).

We added  (c1) to get a self-contained proof for (\ref{wnew}),
and (c3) to draw the attention to cases different from that
shown in Figure 1.
\end{rem}

\begin{thm}
If a tetrahedron has a face which is an equilateral triangle and an 
insphere which touches this face at its center of rotational symmetry, 
then the gap $G = R^2- d^2-3 r^2 - 2 r R$
of the Grace-Danielsson inequality is always zero. \end{thm}

{\bf Proof.} We derive two relations, involving $(w_3, r)$ and
$(w_3, R),$ resp. Like in the proof of Theorem 1, we use lower case letters for vertices, and capitals for the tangent points (e.g. $Z$
is the tangent point of the insphere on the face opposite to $z$).
Denote by $l$ the edge length of the basic equilateral triangle,
i.e. let $l=\|x-y\|=\|y-z\|=\|z-x\|,$ then we have
$ \|(x+y)/2-c\|=l \sqrt{3}/6$ and $\|c-z\|=l\sqrt{3}/3, $
where $c=(0,0,0)$ is the origin. \smallskip

(c1)
Let inc $=(0,0,r)$ be the center of the inscribed sphere,
then $\|w-Z\|$ can be determined from the rectangular triangle 
$\Delta(w,Z,\mathrm{inc})$ using Pythagoras' theorem, 
cf. Figure 1:
\[ \|w-Z\|^2 = \|w-\mathrm{inc}\|^2 - 
\|z-\mathrm{inc}\|^2=(w_3-r)^2-r^2 = w_3\, (w_3-2r), \]
while the similarity of this triangle to
$\Delta(w,\frac{x+y}{2},c)$ implies
\[ \frac{\sqrt{w_3 (w_3-2r)}}{r} = \frac{w_3}{l \sqrt{3}/6}. \]
This immediately gives
\begin{equation} \label{rel1} w_3-2r = \frac{12r^2}{l^2} w_3,
\end{equation}
which implies 
\begin{equation} \label{wnew} w_3=\frac{2l^2r}{l^2-12r^2}. 
\end{equation}

(c2)
Calculating the circumradius by Pythagoras' theorem
applied to the rectangular triangle $\Delta(o, c, z)$ 
gives (see Figure 1):
\[R^2=\|o-c\|^2+\|c-z\|^2=(w_3-R)^2+\frac{l^2}{3} \]
with
\begin{equation} \label{rel2} w_3\, (2R-w_3)=\frac{l^2}{3}. 
\end{equation}
Therefore, by virtue of (\ref{rel1}) and (\ref{rel2}) it follows that
\begin{eqnarray} G&=&(R-r)^2-d^2-4r^2=(R-r)^2-(R+r-w_3)^2
\nonumber \\
 &=&(w_3-2r)\, (2R-w_3)-4r^2 = 4r^2-4r^2 = 0. \nonumber
 \end{eqnarray}
 
(c3)
Note finally, that the order of points $w, o, inc, c$ is
not necessarily that given in Figure 1, hence the relationship
between $w_3, R, r$ and $d$ varies, as well. The precise
formula for the distance $d$ is
\[ d=\begin{cases} R+r-w_3, \quad \mathrm{if} \ \ 0<r<r_{reg}, 
\\ w_3-R-r,  \quad \mathrm{if} \ \ r_{reg}<r<r_{crit} 
\end{cases} \]
where $r_{reg}=\frac{l}{2\sqrt{6}}$ is the inradius of the 
regular tetrahedron (in which case $w_3=R+r$ and $d=0$ hold), 
hence $|d|=|R+r-w_3|,$ and the unified formula 
$d^2=(R+r-w_3)^2$ is valid.
$\square$  \smallskip

In what follows, we examine the analogous  planar problem for triangles. 
Our aim is to obtain a formula for the critical value of the inradius.
\medskip

\noindent {\bf Problem.} Given the line segment $I=[0,1]$
with an interior point $p, \,0<p<1$, 
find the supremum $r_{crit}$ of positive numbers $r,$ for which 
$r$ is the inradius of some triangle with one side equal to  $I.$  
First we illustrate the situation. 

 \begin{figure}[h!]
 \includegraphics[scale=0.7]{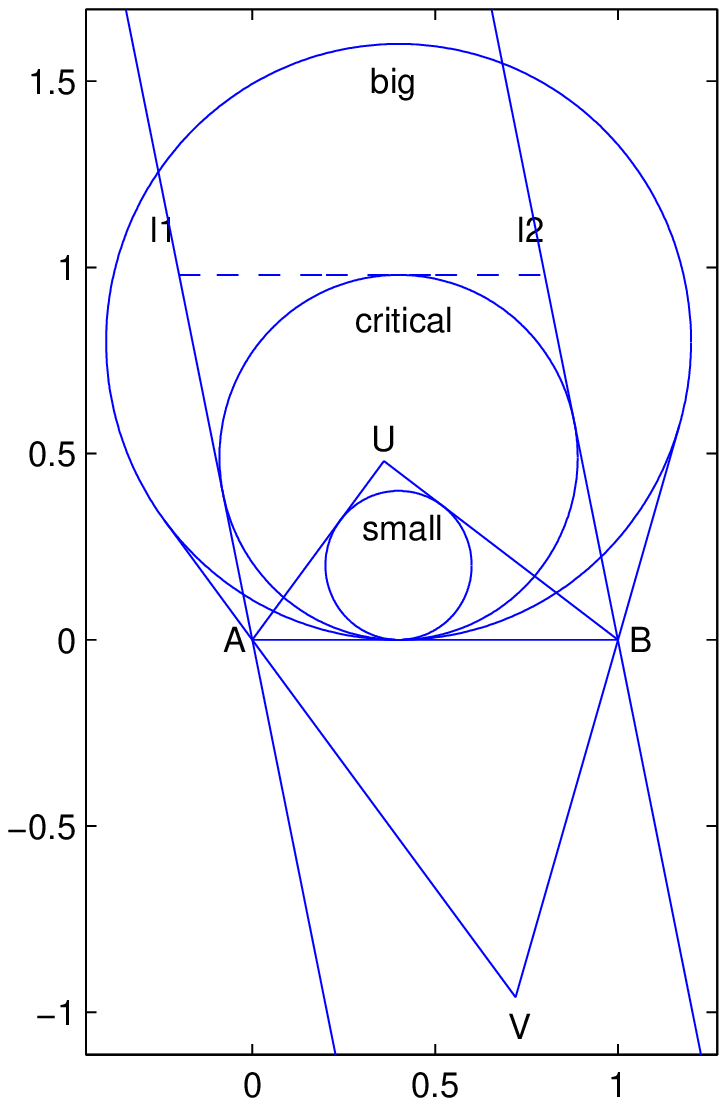}
 \centering
 \centerline{Figure 2. The three cases, p=0.4}
 \end{figure}
 
\begin{rem} Figure 1 below shows a small incircle, resulting in 
triangle $\Delta ABU,$ a critical circle (giving two parallel 
straight lines $l1$ and $l2$ instead of a triangle), 
and a (too) big circle, for which the tangent lines intersect 
at $V,$ on the other (lower) side of the horizontal axis. 
The big circle is then an ex-circle for triangle $\Delta AVB.$
The data  for this plot  are 
\[p=0.4, \ r_{small}=0.2, 
\ r_{crit}=\sqrt{0.24}\approx 0.49, \ r_{big}=0.8.\]
Note that the center $K=(0.4, 0.49)$ of the critical circle is 
quite close to -- but not identical with -- 
the vertex $U=(0.36,0.48)$ of the small right triangle.
\end{rem}

\begin{lem}\label{planar}  For the above planar problem we 
have  $r_{crit} = \sqrt{p(1-p)}.$ \end{lem}
{\bf Proof.}
Triangle $\Delta ABK$ with 
$A=(0,0), \,B=(1,0),\, K=(p,r_{crit})$ is a right triangle. 
To this draw the tangent line to the critical circle, parallel to $AB.$
Then $K$ is the centre of the rhomb bordered by the lines 
$l1, l2$ and the two horizontal tangent lines, hence
$\angle BKA$ is a right angle indeed.

 Using now the well known property:  "the altitude to the hypotenuse 
 is the geometric mean of the two segments of the hypotenuse"
of rectangular triangles, the statement follows.  $\square$ \medskip

After this evasion we go back to three dimensions. In the next
example we calculate the critical inradius, however, in contrast
with the two dimensional case, we can do it only for special data.

\begin{exmp} \label{crit} Let the vertices of the basic triangle be
\[ x= (-\sqrt{2}, -1, 0), \ y=(\sqrt{2}, -1, 0), \ z=(0,1,0), \]
and let the origo be the given interior point.
We show that $r_{crit}=1/\sqrt{2}.$ 
Take for this the sphere $S$ of radius $r=1/\sqrt{2}$ 
centered at $(0, 0, 1/\sqrt{2}),$ and determine the tangent points 
$X, Y, Z$ of the three non-horizontal faces. 
They are
\[  X=\Big( \frac{2\sqrt{2}}{5}, \frac{2}{5},  \frac{2\sqrt{2}}{5} 
\Big),
  \ Y=\Big( -\frac{2\sqrt{2}}{5}, \frac{2}{5},  \frac{2\sqrt{2}}{5} \Big), 
  \ Z=\Big(0, -\frac{2}{3}, \frac{2\sqrt{2}}{3} \Big). \]
The pairwise intersections of the tangent planes give the rays
\[ \Big( \sqrt{2}, b, 2\sqrt{2}\,(1+b) \Big), \ 
   \Big( -\sqrt{2}, b, 2\sqrt{2}\,(1+b) \Big), \
   \Big(0, b, 2\sqrt{2}\,(b-1) \Big) \]
with a free parameter $b.$ Since they share the common direction
 $(0, 1, 2\sqrt{2}),$ the result follows. 

According to Maple, the gap in (\ref{th}) for these 
vertices $(x, y, z)$ is $r^2(1-2r^2),$ showing 
 another evidence for equality $r_{crit}=1/\sqrt{2}.$ 
However the quickest way is 
to show that the points $X, Y, Z$ of tangency with the 
centre $K=(0, 0, 1/\sqrt{2})$ of the sphere are coplanar
(cf.  \cite{coplan}), i.e.
\[ \begin{vmatrix} \frac{2\sqrt{2}}{5} & \frac{2}{5} &       
\frac{2\sqrt{2}}{5} & 1 \cr
  -\frac{2\sqrt{2}}{5} & \frac{2}{5} & \frac{2\sqrt{2}}{5} & 1\cr
0 & -\frac{2}{3} & \frac{2\sqrt{2}}{3} & 1 \cr
0 & 0 & \frac{1}{\sqrt{2}} & 1  \end{vmatrix} = 0. \]
  \end{exmp}

Finally we mention  Pech's method \cite{Pech} proving
Euler's inequality $R\ge 2r$ for triangles (a consequence of 
(\ref{Eul})), to show another idea making use of a computer.
He writes down the known  equations 
 \[r-\frac{2 K}{a+b+c}=0, \quad R-\frac{a b c}{4 K}=0, \quad
  R-2r-k=0,\]
as well as Heron's formula
\[16 K^2-(a+b+c)(a+b-c)(a-b+c)(-a+b+c)=0, \]
including $r, \ R,$ the area $K,$ the three sides $a, b, c$ of a 
triangle, and a slack variable $k.$
Using the CoCoA (short for Computations in Commutative 
Algebra) system he finds that $R\ge 2 r$ holds iff
\[a^3-a^2 b-a b^2 +b^3 -a^2 c +3 abc -b^2 c-a c^2 -b c^2 
+c^3 \ge 0,\]
which is easily shown by observing that this polynomial equals
\[\tfrac{1}{2}\, [(a + b - c) (a - b)^2 +( b + c - a)(b - c)^2 
+( c + a - b)(c - a)^2]. \]

Note that Pech's method is much simpler than ours, 
thanks to its coordinate-free approach, 
however, it applies to the planar case $n=2$, and
-- on the other hand --, it
does not concern the distance $d.$ 
For $n=3$ it would be a challenge to express  $d$  
by help of lengths only -- without using coordinates. 

\section*{Acknowledgement}
My special thanks to the reviewer for his helpful and constructive 
comments and suggestions to improve the paper.

\appendix

\section{Appendix}
The polynomials $v_1$ and $v_2$ have the following, 
fairly symmetrical form:
\begin{eqnarray*}
v_1&=& c_1\, (c_1^2+c_2^2)\,a0-2\,c_1^2\, a0\, o_1 \\
 &+& y_1z_1(y_2-z_2)(c_1^2-x_1^2-(c_2-x_2)^2) \\
 &+& z_1x_1(z_2-x_2)(c_1^2-y_1^2-(c_2-y_2)^2) \\
 &+& x_1y_1(x_2-y_2)(c_1^2-z_1^2-(c_2-z_2)^2) \\
 &+&  c_1^2(x_2^2(z_2-y_2)+y_2^2(x_2-z_2)+z_2^2(y_2-x_2)) \\
 &+&  c_2^2(x_1^2(z_2-y_2)+y_1^2(x_2-z_2)+z_1^2(y_2-x_2)) \\
 &+&  x_1^2(y_2-z_2)(c_1(y_1+z_1)+c_2(y_2+z_2)-y_2z_2) \\
 &+&  y_1^2(z_2-x_2)(c_1(z_1+x_1)+c_2(z_2+x_2)-z_2x_2) \\
 &+&  z_1^2(x_2-y_2)(c_1(x_1+y_1)+c_2(x_2+y_2)-x_2y_2) \\
 &+&  2c_1c_2(x_1x_2(y_2-z_2)+y_1y_2(z_2-x_2)+z_1z_2(x_2-y_2)) \\
 &+&  2c_1c_2(x_1(z_2^2-y_2^2)+y_1(x_2^2-z_2^2)+
         z_1(y_2^2-x_2^2)) \\
 &+&  c_1(x_1x_2(z_2^2-y_2^2)+y_1y_2(x_2^2-z_2^2)+
         z_1z_2(y_2^2-x_2^2)) \\
 &+& 3c_1x_1y_2z_2(y_2-z_2)+x_2y_1z_2(z_2-x_2)+
        x_2y_2z_1(x_2-y_2)), \end{eqnarray*}
and
\begin{eqnarray*}
v_2 &=& c_2\,(c_1^2+c_2^2)\,a0-2\,c_2^2\,a0\,o_2 \\
 &=&  y_2z_2(y_1-z_1)((c_1-x_1)^2-c_2^2+x_2^2) \\
 &=&  x_2z_2(z_1-x_1)((c_1-y_1)^2-c_2^2+y_2^2) \\
 &=&  x_2y_2(x_1-y_1)((c_1-z_1)^2-c_2^2+z_2^2) \\
 &=&  c_1^2(x_2^2(y_1-z_1)+y_2^2(z_1-x_1)+z_2^2(x_1-y_1)) \\
 &=&  c_2^2(x_1^2(y_1-z_1)+y_1^2(z_1-x_1)+z_1^2(x_1-y_1)) \\
 &=&  x_2^2(z_1-y_1)(c_1(y_1+z_1)+c_2(y_2+z_2)-y_1z_1) \\
 &=&  y_2^2(x_1-z_1)(c_1(z_1+x_1)+c_2(z_2+x_2)-z_1x_1) \\
 &=&  z_2^2(y_1-x_1)(c_1(x_1+y_1)+c_2(x_2+y_2)-x_1y_1) \\
 &=&  2c_1c_2(x_1x_2(z_1-y_1)+y_1y_2(x_1-z_1)+z_1z_2(y_1-x_1)) \\
 &=&  2c_1c_2(x_2(y_1^2-z_1^2)+y_2(z_1^2-x_1^2)+z_2(x_1^2-y_1^2)) \\
 &=&  c_2(x_1x_2(y_1^2-z_1^2)+y_1y_2(z_1^2-x_1^2)+
         z_1z_2(x_1^2-y_1^2)) \\
 &=&  3c_2(x_1^2(y_2z_1-y_1z_2)+y_1^2(x_1z_2-x_2z_1)+
         z_1^2(x_2y_1-x_1y_2)). \end{eqnarray*}

\end{document}